\documentclass[12pt]{article}
\usepackage{amsmath}\usepackage{amssymb}

\date{}

\title{Two examples concerning almost continuous functions}

\author{
Krzysztof Ciesielski%
\thanks{AMS classification numbers: Primary  26A15, 26A30;
        Secondary  03E50. \endgraf
  Key words and phrases: additive, almost continuous,
         extendability, and SCIVP functions. \endgraf
  The first author was partially supported by NSF Cooperative
         Research Grant INT-9600548 with its Polish part financed by KBN.}\\
{\footnotesize
    Department of Mathematics, West Virginia University,}\\
{\footnotesize  Morgantown, WV 26506-6310} \\
{\footnotesize  KCies@wvnvms.wvnet.edu}
\and
Andrzej Ros{\l}anowski\\
{\footnotesize Department of Mathematics and Computer Science}\\
{\footnotesize 1910 University Dr}\\
{\footnotesize Boise State University}\\
{\footnotesize Boise, ID 83725}\\
{\footnotesize roslanow@math.idbsu.edu}
}

\newtheorem{theorem}{Theorem}[section]
\newtheorem{corollary}[theorem]{Corollary}
\newtheorem{proposition}[theorem]{Proposition}
\newtheorem{lemma}[theorem]{Lemma}
\newtheorem{problem}{Problem}[section]
\newtheorem{example}{Example}[section]
\newtheorem{definition}{Definition}
\newtheorem{conjecture}{Conjecture}

\newcommand{\thm}[2]{\begin{theorem}\label{#1}#2\end{theorem}}
\newcommand{\cor}[2]{\begin{corollary}\label{#1}#2\end{corollary}}
\newcommand{\prop}[2]{\begin{proposition}\label{#1}#2\end{proposition}}
\newcommand{\lem}[2]{\begin{lemma}\label{#1}#2\end{lemma}}
\newcommand{\pr}[2]{\begin{problem}\label{#1}#2\end{problem}}

\newcommand{\T}{{\cal T}}

\newcommand{\real}{{\mathbb R}}
\newcommand{\bR}{\real}
\newcommand{\rational}{{\mathbb Q}}
\newcommand{\bQ}{\rational}
\newcommand{\integer}{{\mathbb Z}}
\newcommand{\bZ}{\integer}

\def\continuum{{\mathfrak c}}
\def\co{\continuum}
\def\con{\continuum}

\newcommand{\op}{\operatorname}
\newcommand{\Add}{{\op{Add}}}
\newcommand{\AC}{{\op{AC}}}
\newcommand{\CIVP}{{\op{CIVP}}}
\newcommand{\SCIVP}{{\op{SCIVP}}}
\newcommand{\ext}{{\op{Ext}}}
\newcommand{\Conn}{{\op{Conn}}}

\def\lin{{\operatorname{LIN}_\rational}}

\newcommand{\cl}{{\op{cl}}}
\newcommand{\Fr}{{\op{bd}}}
\newcommand{\proj}{{\op{proj}}}
\newcommand{\rng}{{\op{rng}}}
\newcommand{\diam}{{\op{diam}}}

\newcommand{\bd}{{\op{bd}}}

\def\proof{{\sc \noindent Proof }}
\newcommand{\Proof}{\proof}
\def\qed{\hfill\vrule width 6pt height 6pt depth 0pt\vspace{0.1in}}
\newcommand{\QED}{\qed}
\def\la{\langle}
\def\ra{\rangle}

\newcommand{\baire}{{\omega^{\textstyle\omega}}}
\newcommand{\rest}{\restriction}

\begin{document}

\maketitle

\begin{abstract}

In this note we will construct, under the assumption that
union of less than continuum many meager subsets of $\real$ is meager in
$\real$,
an additive
connectivity function $f\colon\real\to\real$
with Cantor intermediate value property
which is not almost
continuous. This gives a partial answer to a question of
D.~Banaszewski~\cite{Ban}.
(See also \cite[Question~5.5]{GN}.)
We will also show that every extendable function $g\colon\real\to\real$
with a dense graph satisfies the following  stronger version of the SCIVP
property: for every $a<b$ and every perfect set $K$ between
$g(a)$ and $g(b)$ there is a perfect set $C\subset(a,b)$ such that
$g[C]\subset K$ and $g\restriction C$ is continuous {\em strictly increasing}.
This property is used to construct a ZFC example of an
additive almost
continuous function $f\colon\real\to\real$ which has the strong Cantor
intermediate value property but is not extendable.
This answers a question of H.~Rosen~\cite{Ro}. This also
generalizes Rosen's result~\cite{Ro} that
a similar (but not additive)
function exists under the assumption of
the continuum hypothesis, and gives a full answer to \cite[Question~3.11]{GN}.
\end{abstract}

\section{Preliminaries}\label{sec0}

Our terminology is standard and follows \cite{CiBook}.  We consider only
real-valued functions of one or two real variables. No distinction is made
between a
function and its graph. By $\real$ and $\rational$ we denote the set of all
real and rational numbers, respectively.  We will consider $\real$ and
$\real^2$ as linear spaces over $\rational$.  In particular, for a subset $X$
of either $\real$ or $\real^2$ we will use the symbol $\lin(X)$ to denote the
smallest linear subspace (of $\real$ or $\real^2$) over $\rational$ that
contains $X$.  Recall also that if $D\subset\real$ is linearly independent
over $\rational$ and $f\colon D\to\real$ then $F=\lin(f)\subset\real^2$ is an
additive function (see definition below) from $\lin(D)$ into $\real$.  Any
linear basis of $\real$ over $\rational$ will be referred as a {\em Hamel
basis}. By a Cantor set we mean any nonempty perfect nowhere dense subset of
$\real$.

The ordinal numbers will be identified with the sets of all their
predecessors and cardinals with the initial ordinals. In particular
$2=\{0,1\}$, and the first infinite ordinal $\omega$ number is equal to the
set of all natural numbers $\{0,1,2,\ldots\}$.  The family of all functions
from a set $X$ into $Y$ is denoted by $Y^X$.  The symbol $|X|$ stands for the
cardinality of a set $X$. The cardinality of $\real$ is denoted by $\co$ and
referred as {\em continuum}.  A set $S\subset\real$ is said to be {\em
$\co$-dense\/} if $|S\cap(a,b)|=\co$ for every $a<b$.  The closure of a set
$A\subseteq\bR$ is denoted by $\cl(A)$, its boundary
by $\Fr(A)$,
and its diameter by $\diam(A)$.  For a set $A\subseteq X\times Y$ and points
$x\in X$ and $y\in Y$ we let $(A)_x = \{y\in Y\colon\la x,y\ra\in A\}$ and
$(A)^y = \{x\in X\colon\la x,y\ra\in A\}$.  In a similar manner we define
$(A)_{\la x,y\ra}$ and $(A)^z$ for a set $A\subseteq X\times Y\times Z$.

We will use also the following terminology~\cite{GN}. A function $f\colon
\real\to\real$
\begin{itemize}
\item is {\em additive\/}
if $f(x+y)=f(x)+f(y)$ for every $x,y\in\real$;

\item is {\em almost continuous\/} (in sense of Stallings)
if each open subset of $\real\times\real$ containing the graph of
$f$ contains also a continuous function from $\real$ to $\real$~\cite{Stal};

\item has the {\it Cantor intermediate value property\/} if for every
$x,y\in\real$ and for each Cantor set $K$ between $f(x)$ and $f(y)$ there is
a Cantor set $C$ between $x$ and $y$ such that $f[C]\subset K$;

\item has the {\it strong Cantor intermediate value property\/} if for every
$x,y\in\real$ and for each Cantor set $K$ between $f(x)$ and $f(y)$ there is
a Cantor set $C$ between $x$ and $y$ such that $f[C]\subset K$ and
the restriction $f\restriction C$ of $f$ to $C$ is continuous;

\item is an {\it extendability\/} function provided
there exists  a connectivity function $F\colon\real\times[0,1]\to\real$
such that $f(x)=F(x,0)$ for every $x\in\real$, where

\item for a topological space $X$ a function $f\colon X\to\real$ is a
{\it connectivity\/} function if the graph of the restriction  $f\restriction
Z$ of $f$ to $Z$ is connected in $Z\times\real$ for any connected subset $Z$
of~$X$.
\end{itemize}
The above classes of functions (from $\real$ to $\real$) will be denoted by
$\Add$,
$\AC$, $\CIVP$, $\SCIVP$, $\ext$, and $\Conn$, respectively.

Recall that if the graph of $f\colon\real\to\real$ intersects every closed
subset $B$ of $\real^2$ which projection $\proj(B)$ onto the $x$-axis has
nonempty interior then $f$ is almost continuous.  (See e.g.~\cite{N1}.)
Similarly, if the graph of $f\colon\real\to\real$ intersects every compact
connected subset $K$ of $\bR^2$ with $|\proj(K)|>1$ then $f$ is connectivity.

We will finish this section with the following well known fact.
(See \cite[Theorem 4.A.12]{ChVo}, \cite[\& 47III]{Kur}, or
\cite[Ch. V, sec. 2]{Nadler}.)

\prop{prop:BBT}{{\rm (Boundary Bumping Theorem)}
If $U$ is a nonempty open proper
subset of a compact connected Hausdorff space $K$ and $C$ is
a connected component of $U$ then $\cl_K(C)\cap\Fr_K(U)\neq \emptyset$.
In particular every connected component of $U$ has more than one point.
}

\section{Additive connectivity function on $\bR$ which is not almost continuous}

We start this section with recalling the following construction of
Roberts~\cite{Ro65}
of zero-dimensional closed subset $Z_0$ of $[0,1]^2$ which is intersected
by a graph of every continuous function $f\colon[0,1]\to[0,1]$.
Let $C\subset[0,1]$ be a Cantor set of Lebesgue measure $1/2$.
(Roberts defines it as $C=\bigcap_{n<\omega}C_n$, where $C_0=[0,1]$,
each $C_n$ is the union of $2^n$ disjoint intervals, and
$C_{n+1}$ is obtained from $C_n$ by taking out of each of these $2^n$ intervals
a concentric open interval of length $1/2^{2n+2}$.)
Define $x,y\colon[0,1]\to[0,1]$ by $x(t)=2m(C\cap[0,t])$,
where $m$ is a Lebesgue measure, and
$y(t)=4 m(C\cap[0,t])-t=2 x(t)-t$.
Then $F_0\colon[0,1]\to[0,1]^2$, $F_0(t)=\la x(t),y(t)\ra$, is a continuous
embedding, so $M_0=F_0[[0,1]]$ is an arc joining $\la 0,0\ra$ with $\la 1,1\ra$.
Note that each component interval $I$ of $[0,1]\setminus C$ is mapped by $F_0$
onto an open vertical segment $F_0[I]$. The set $Z_0$ defined as $F_0[C]$.
It is equal to the arc $M_0$ from which all vertical segments $F_0(I)$ are
removed.
Note also that an arc $F_0[I]$ has been removed from the section
$(M_0)_x$ if and only if $x\in D_0$, where $D_0$ is the set of all
dyadic numbers ($x=k/2^n$) from $(0,1)$.
Moreover, $|(Z_0)_x|=2$ for $x\in D_0$ and $(Z_0)_x=(M_0)_x$ is a singleton
for all other $x$ from $[0,1]$.

For what follows we will need the following version of
this construction, where $\bar C=\integer+C$.

\lem{roberts}{
Let $X$ be a countable dense subset of $(-1,1)$.
Then there exists an embedding
$F=\la F_0,F_1\ra\colon\real\to(-1,1)\times\real$ such that
$F_0$ is non-decreasing,
\begin{description}
\item[(a)] an open arc $M=F[\real]$ is closed in $\real^2$,
\item[(b)] if $Z=F[\bar C]\subset M$ then
           $g\cap Z\neq\emptyset$ for every continuous
           $g:[-1,1]\to\bR$,
\item[(c)] $Z_x=M_x$ is a singleton for all $x\in(-1,1)\setminus X$, and
\item[(d)] for each  $x\in X$ the section $M_x$ is a non-trivial
           closed interval and
           $Z_x$ consists of the two endpoints of that interval.
\end{description}
}

\Proof Let $F_0$ be Roberts' function defined above.
Define $F_1\colon\real\to\real^2$ by putting $F_1(n+x)=\la n,n\ra+F_0(x)$
for every $n\in\integer$ and $x\in[0,1)$. Then $F_1$ is a continuous embedding
extending $F_0$. Also choose an order isomorphism $h\colon\real\to(-1,1)$
such that $h[\integer+D_0]=X$ and define a homeomorphism
$H\colon\real^2\to(-1,1)\times\real$ by $H(x,y)=\la h(x),y\ra$.
It easily follows from the properties of $F_0$
that $F=H\circ F_1$ satisfies (a)-(d).
\QED

Note that by (b) of Lemma~\ref{roberts} if the graph of
$f\colon\real\to\real$ is disjoint with $Z$ then
$f$ is not almost continuous, since then the set $U=\real^2\setminus Z$
is an open set containing $f$ which does not
contain any continuous function $g\colon\real\to\real$.
Thus the main idea of the next theorem is to construct an additive
connectivity function with the graph disjoint with $Z$.

In our argument it will be also convenient to use the following
easy lemma.

\lem{lem:Hamel}{Let $\{I_\alpha\colon\alpha<\co\}$
be an enumeration, with possible repetitions, of all
nonempty open intervals in $\real$. Then there exists
a family of pairwise disjoint perfect sets
$\{P_\alpha\subset I_\alpha\colon\alpha<\co\}$ such that
$P=\bigcup_{\alpha<\co}P_\alpha$ is meager in $\real$ and
linearly independent over $\rational$.
Moreover, we can assume that there is a
meager $F_\sigma$-set $S$ containing $P$ such that
$S=\lin(S)$ and $S$ is of co-dimension continuum.
}

\proof Take a linearly independent perfect subset $K$ of $\real$.
(See e.g. \cite[thm.~2, Ch.~XI sec.~7]{Kucz}.)
Partition $K$ into perfect sets $\{F,H,L\}$ and further partition
$F$ into pairwise disjoint perfect sets $\{F'_\alpha\colon\alpha<\co\}$.
Choose a countable subset $H_0=\{x_n\colon n<\omega\}$ of $H$
and for every $\alpha<\co$ choose a sequence
of non-zero rational numbers $\la q^\alpha_n\colon n<\omega\ra$
such that $F_\alpha=\bigcup_{n<\omega}q^\alpha_n\cdot x_n + F'_\alpha$
is dense in $\real$. Then the sets $F_\alpha$ are pairwise disjoint
and $\bigcup_{\alpha<\co}F_\alpha$ is linearly
independent over $\rational$.
For every $\alpha<\co$ choose perfect $P_\alpha\subset F_\alpha\cap I_\alpha$.
Then $P_\alpha$'s are pairwise disjoint and
$P=\bigcup_{\alpha<\co}P_\alpha\subset\rational\cdot H_0 + F$
is meager. Also if
$S=\lin(H\cup F)$ then $S=\lin(S)$, and it is an $F_\sigma$-set.
It is of co-dimension continuum (so meager)
since $S$ is disjoint with $L$. \qed

\thm{AddAndConnNotAC}{
If union of less than $\con$ many meager subsets of $\real$ is meager in
$\bR$ then there exists an $f\in\Add\cap\CIVP\cap\Conn\setminus\AC$.
}

\Proof Let $\la\la I_\alpha,C_\alpha\ra\colon\alpha<\co\ra$ be a list of all
pairs $\la I,C\ra$ such that $I$ is a nonempty open interval in $\real$ and
$C$ is a perfect subset of $\real$ and take
$\{P_\alpha\subset I_\alpha\colon\alpha<\co\}$
as in Lemma~\ref{lem:Hamel}.

Let $\{C,D\}$ be a partition of $\co\setminus\omega$
onto sets of cardinality continuum.
Take an enumeration
$\{K_\xi\colon\xi\in D\}$ of the family of all compact
connected subsets $K$ of $\bR^2$ with $|\proj[K]|=\con$.
Also, let $H$ be a Hamel basis containing
$P=\bigcup_{\alpha<\co}P_\alpha$ such that
there is a countable set $X\subset (H\setminus P)\cap(-1,1)$ dense
in $(-1,1)$.
Let $Z$ be as in Lemma~\ref{roberts} for this $X$ and
$\{h_\xi\colon\xi\in C\}$ be an enumeration of $H$.
By induction on $\xi<\con$ we will
choose functions $f_\xi$ from finite subsets $H_\xi$ of $H$ into $\real$
such that for every $\xi<\con$ the
following conditions hold.
\begin{enumerate}
\item[(i)]   $H_\xi\cap\bigcup_{\zeta<\xi}H_\zeta=\emptyset$.

\item[(ii)] If $\xi\in C$ then $h_\xi\in \bigcup_{\zeta\leq\xi}H_\zeta$.

\item[(iii)] If $\xi\in D$ then
$K_\xi\cap\lin(\bigcup_{\zeta\leq\xi}f_\zeta)\neq\emptyset$.

\item[(iv)] $Z\cap\lin(\bigcup_{\zeta\leq\xi}f_\zeta)=\emptyset$.

\item[(v)] If
$x\in H_\xi\cap P_\alpha$ for some $\alpha<\co$ then $f_\xi(x)\in C_\alpha$.
\end{enumerate}

Before we describe the inductive construction note first how it
can be used to construct a function as desired.
First notice that, by (i) and (ii), $\bigcup_{\xi<\co}f_\xi$
is a function from $H$ into $\real$. Thus
$f=\lin(\bigcup_{\xi<\co}f_\xi)$ is an additive function from $\real$ to
$\real$.
It is connectivity by (iii). It is not almost continuous
by (iv) and remark after Lemma~\ref{roberts}.
It has Cantor intermediate value property by (v) and the choice of
$\la I_\alpha,C_\alpha\ra$.

The main difficulty in our inductive construction will be the preservation of
condition (iv). To handle this easier note that if
$g$ is an additive function from $E\subset\real\setminus\{x\}$
into $\real$ such that
$Z\cap g=\emptyset$ then $Z\cap\lin(g\cup\{\la x,y\ra\})=\emptyset$
if and only if
\begin{equation}\label{eq:1}
\la x,y\ra\notin \bigcup\{q Z+\la p,g(p)\ra\colon p\in E\ \&\ q\in\rational\}.
\end{equation}
In particular, if $x$ is fixed, than
$Z\cap\lin(g\cup\{\la x,y\ra\})=\emptyset$
if and only if
\begin{equation}\label{eq:2}
y\notin \bigcup\{(q Z+\la p,g(p)\ra)_x\colon p\in E\ \&\ q\in\rational\}.
\end{equation}

We will make the construction in two main steps.
First we will construct the functions $f_n$ for $n<\omega$.
For this choose an enumeration $\{x_n\colon n<\omega\}$ of $X$.
We put $H_n=\{x_n\}$ and define $f_n(x_n)$
inductively such that
\begin{equation}\label{eq:vert}
\la x_n,f_n(x_n)\ra\in M\setminus Z,
\end{equation}
where $M$ is the set from Lemma~\ref{roberts}.

To see that such a choice can be made, note first that (i) is satisfied,
and (ii), (iii), and (v) are satisfied in void.
Thus, we have to take care only of the condition (iv).
However, for each $n<\omega$ we have an entire interval
of possible choices for $f_n(x_n)$ (see Lemma~\ref{roberts}(d))
while, by (\ref{eq:2}),
there is only a countable many exceptional points we have to avoid.
(Since $|Z_{x_n}|=2$ and $E=\lin(\{x_i\colon i<n\})$ in this case.)

Now, assume that for some infinite $\xi<\co$ the sequence
$\la f_\zeta\colon \zeta<\xi\ra$ has been already constructed.
Put $g=\lin(\bigcup_{\zeta<\xi}f_\zeta)$ and let $E$ be its domain.

First consider case when $\xi\in C$. If $h_\xi\in\bigcup_{\zeta<\xi}H_\zeta$
we put $f_\xi=H_\xi=\emptyset$.
So, assume that $h_\xi\notin\bigcup_{\zeta<\xi}H_\zeta$
and put $H_\xi=\{h_\xi\}$.
If $h_\xi\in P_\alpha$ for some $\alpha<\co$ put $P=C_\alpha$.
Otherwise put $P=\real$. Then (i) and (ii) are satisfied
and (v) will hold if we choose $f_\xi(h_\xi)\in P$.
To have (iv) by (\ref{eq:2}) it is enough to choose
$f_\xi(h_\xi)$ from outside of a set
$\bigcup\{(q Z+\la p,g(p)\ra)_{x_\xi}\colon p\in E\ \&\ q\in\rational\}$,
which has cardinality less than continuum.

So, assume that $\xi\in D$. Let $S$ be as in Lemma~\ref{lem:Hamel}
and put
$T_0=\lin(S\cup\bigcup_{\zeta<\xi}H_\zeta)$. Then
$T_0\neq\real$ since $\lin(S)$ is of co-dimension
continuum. Moreover $T_0$ is a union of less than continuum many
meager sets $\lin(S\cup A)$, where $A$ is a finite subset of
$\bigcup_{\zeta<\xi}H_\zeta$.
Thus, by our assumption, $T_0$ is meager. Let $T$ be a meager $F_\sigma$-set
containing $T_0$.
Our next main objective will be to show that
either we already have $K_\xi\cap g\neq\emptyset$ or
we can find
\begin{equation}\label{eq:main}
\la x,y\ra\in K_\xi\setminus\left(
(T\times\real)\cup\bigcup\{q Z+\la p,g(p)\ra\colon p\in E\ \&\ q\in\rational\}
\right).
\end{equation}

Before we argue for it, first note how this will finish the construction.
If $K_\xi\cap g\neq\emptyset$ we can put $f_\xi=H_\xi=\emptyset$.
So, assume that we can find $\la x,y\ra$ as in (\ref{eq:main}).
Take a minimal subset
$\{k_0,\ldots,k_m\}$ of $H\setminus\bigcup_{\zeta<\xi}H_\zeta$
such that $x\in\lin(\{k_0,\ldots,k_m\}\cup\bigcup_{\zeta<\xi}H_\zeta)$.
We will define $f_\xi$ on $H_\xi=\{k_0,\ldots,k_m\}$ such that
$\la x,y\ra\in\lin(g\cup f_\xi)$, implying (iii),
while preserving (iv) and (v). First, for $i\leq m$ let $P^i$ be equal to
$C_\alpha$ if $k_i\in P_\alpha$ for some $\alpha<\co$
and equal to $\real$ otherwise. To preserve (v) we have to choose
$f_\xi(k_i)\in P^i$.
Next note that
$H_\xi\not\subset P$ since $x\notin T\supset
\lin(P\cup\bigcup_{\zeta<\xi}H_\zeta)$.
Assume that $k_m\notin P$. Thus $P^m=\real$.
Note that, by (\ref{eq:1}),
$\bar g=\lin(g\cup\{\la x,y\ra\})$ is disjoint with $Z$.
Proceeding as in case when $\xi\in C$ and
using (\ref{eq:2}) we can inductively choose for every $i<m$
a value $f_\xi(k_i)\in P^i$ such that
$h=\lin(\bar g\cup\{\la k_i,f_\xi(k_i)\ra \colon i<m\})$
is disjoint with $Z$. Then function $h$ is already defined on $k_m$
and we can put $f_\xi(k_m)=h(k_m)\in\real=P^m$.
Clearly such $f_\xi$ satisfies (iv) and (v).

To argue for (\ref{eq:main}) we will consider three cases.

\medskip

\noindent{\sc Case 1:}
$\emptyset\neq (I\times\bR)\cap(qM+v)\subset K_\xi$
for some $v=\la v_0,v_1\ra\in g$, $q\in \bQ\setminus\{0\}$,
and an open interval $I$.
Then $K_\xi\cap g\neq\emptyset$.

\medskip

Indeed $\frac{1}{q}(I-v_0)$ is an open interval intersecting $(-1,1)$
and we
find $n<\omega$ such that $x_n\in H_n\cap\frac{1}{q}(I-v_0)$.
By (\ref{eq:vert}) we have $\la x_n,g(x_n)\ra\in M\setminus Z$.
Therefore
$\la q x_n+v_0,g(q x_n+v_0)\ra=
q\la x_n,g(x_n)\ra+v\in (qM+v)\cap(I\times\bR)\subset K_\xi$.

\medskip

\noindent{\sc Case 2:}
There exists an $x\in\{z\in\bR\colon |(K_\xi)_z|=\con\}\setminus T$.

\medskip

Choose $y\in (K_\xi)_x\setminus
\bigcup\{(q Z+\la p,g(p)\ra)_x\colon p\in E\ \&\ q\in\rational\}$.
Then $\la x,y\ra$ satisfies~(\ref{eq:main}).

\medskip

\noindent{\sc Case 3:}
Neither Case 1 nor Case 2 hold.

\medskip

Define $Y$ as $K_\xi\setminus(T\times\bR)$. Then $Y$ is a $G_\delta$
subset of $K_\xi$ so it is a Polish space.
Notice also that, since we are not in Case~2, every vertical section of
$Y$ is at most countable.
We will prove that
\begin{equation}\label{eq:5}
q Z + v \ \mbox{ is meager in $Y$
for every $v\in g$ and $q\in\rational$}.
\end{equation}
This clearly implies the possibility of a choice as in (\ref{eq:main})
since $\real$ (and so, a Polish space $Y$) is not a union
of less than continuum many meager sets.

To prove (\ref{eq:5}) fix $v=\la v_0,v_1\ra\in g$,
$q\in\rational\setminus\{0\}$,
and an open set $U\subset\real^2$ such that $U\cap Y\neq\emptyset$.
We have to show that $U\cap Y\setminus (q Z + v)\neq\emptyset$.
So, fix $p=\la x,y\ra\in U\cap Y$ and an open set $V$ containing $p$ such that
$\cl(V)\subset U$. Let $C_0$ be a connected component
of $K_\xi\cap V$ containing $x$. Then, by Proposition~\ref{prop:BBT},
$C_0$ has more than one point. Consider a compact connected set
$K=\cl(C_0)\subset\cl(V)\subset U$. Then $p\in K$ and $\proj(K)$
is a non-trivial interval, say $[c,d]$, since $K_x\subset(K_\xi)_x$
is at most countable. Thus, it is enough to prove that
$K\setminus\left((T\times\real)\cup(q Z + v)\right)\neq\emptyset$
which follows easily from the following property:
\begin{equation}\label{eq:6}
|\proj(C)|=\co\  \mbox{ for some connected component $C$ of
$K\setminus(q M + v)$},
\end{equation}
where $M$ is an arc from Lemma~\ref{roberts} containing $Z$.

By way of contradiction assume that (\ref{eq:6}) is false. Then
every connected component of $K\setminus(q M + v)$ is vertical.
Note that there exists a number $r\in (q X + v_0)\cap(c,d)$ such that
the vertical section $\{r\}\times (q M+v)_r$ of $q M + v$
is not contained in $K$, since otherwise
we would have
\[
((c,d)\times\real)\cap(q M + v)\subset
\cl\left(\bigcup_{r\in (q X + v_0)\cap(c,d)}\{r\}\times (q M+v)_r\right)\subset
K\subset K_\xi
\]
contradicting the fact that Case~1 does not hold.
Let
$a<b$ be such that $\la r,a\ra$ and $\la r,b\ra$ are the endpoints
of the vertical segment $\{r\}\times[a,b]$ of $q M + v$ above $r$,
i.e., such that $(q M + v)_r=[a,b]$. Since $\{r\}\times[a,b]$
is not a subset of a compact set $K$, we can find $s\in(a,b)$
such that $\la r,s\ra\notin K$. Take an $\varepsilon_0>0$ such that
\begin{enumerate}
\item[$(\alpha)$] $\varepsilon_0<\frac{1}{4}\min\{s-a,b-s,r-c,d-r\}$ and
\item[$(\beta)$]  the closed rectangle $[r-\varepsilon_0,r+\varepsilon_0]\times
[s-\varepsilon_0,s+\varepsilon_0]$ is disjoint from $K$.
\end{enumerate}
It follows from Lemma~\ref{roberts}
(in particular, the fact that $F_0$ is non-decreasing)
that we may find a positive $\varepsilon_1<
\varepsilon_0$ such that either
\begin{enumerate}
\item[$(\gamma)$] $(\forall x\in (r,r+\varepsilon_1])(\forall y\in\big(qM+v
\big)_x)(a-\varepsilon_0<y<a+\varepsilon_0)$, and
\item[$(\delta)$] $(\forall x\in [r-\varepsilon_1,r))(\forall y\in\big(qM+v
\big)_x)(b-\varepsilon_0<y<b+\varepsilon_0)$,
\end{enumerate}
or symmetrical conditions interchanging $(r,r+\varepsilon_1]$,
$[r-\varepsilon_1,r)$ hold. Without loss of generality we may assume that we
have the clauses $(\gamma)$, $(\delta)$ as formulated above.
(For $Z$ and $M$ as constructed in Lemma~\ref{roberts} this happens when
$q>0$.) Consider the set
\[
D_{\varepsilon_1}\stackrel{\rm def}{=}
([r-\varepsilon_1,r+\varepsilon_1]\times\{s\})\cup
\{\la r+\varepsilon_1,y)\ra\colon y\geq s\}\cup
\{(\la r-\varepsilon_1,y\ra\colon y\leq s\}.
\]
We claim that $D_{\varepsilon_1}\cap K=\emptyset$. Why? Suppose that
$\la x,y\ra\in D_{\varepsilon_1}\cap K$.
By the choice of $\varepsilon_0$ (clause $(\beta)$)
we know that either $x=r+\varepsilon_1$ and $y>s$,
or $x=r-\varepsilon_1$ and $y<s$.
The two cases are handled similarly, so suppose that the first one takes place.
By the choice of $\varepsilon_1$ (clause $(\gamma)$) we know that
$\la x,y\ra\notin qM+v$ (as $y>s>a+\varepsilon_0$).
We have assumed that each connected
component of $K\setminus (qM+v)$ is contained in a vertical line, so look at
the connected component $C_{\la x,y\ra}$ of $K\setminus(qM+v)$
to which $\la x,y\ra$ belongs.
By Proposition~\ref{prop:BBT}
we know that $\cl(C_{\la x,y\ra})\cap(qM+v)\neq\emptyset$.
Hence, by clause $(\gamma)$, we conclude that
$\la x,s\ra\in C_{\la x,y\ra}$ (remember
$y>s>a+\varepsilon_0$), a contradiction with clause $(\beta)$.

To obtain a final contradiction note that $D_{\varepsilon_1}$
separates non-empty subsets $K\cap(\{c\}\times\real)$
and $K\cap(\{d\}\times\real)$, which contradicts connectedness of $K$.
The proof is complete. \QED

It is also worth to mention that essentially the same proof as above gives
the following theorem with a slightly weaker set theoretical assumption.

\thm{AddAndConnNotAC2}{If $\real$ is not a union of less than continuum
many of its meager subsets
then there exists an $f\in\Add\cap\Conn\setminus\AC$.
}

{\sc\noindent Sketch of proof } The argument can be obtained by the following
modification of the proof of Theorem~\ref{AddAndConnNotAC}.
Repeat the proof with replacing sets $S$, $P$, and
$P_\alpha$'s with the empty set.
Then (v) is always satisfied in void
and $T$ will become $\lin(\bigcup_{\zeta<\xi}H_\zeta)$,
which has cardinality less than $\co$, but
certainly does not have to be $F_\sigma$.
Then we note that the set $A=\{z\in\bR\colon |(K_\xi)_z|=\con\}$
is analytic, so it is either countable, or has cardinality continuum.
Thus, if Case~2 does not hold then $A$ is countable.
The proof is finished when we replace the set $X$ from the
proof of Theorem~\ref{AddAndConnNotAC}
with $X=K_\xi\setminus(A\times\real)$
and notice that the sets $\{z\}\times \real$ with
$z\in\lin(\bigcup_{\zeta<\xi}H_\zeta)$ are meager in $X$. \qed

We will finish this section with the following open problems.

\pr{pr1}{ Does there exist a ZFC example of an additive connectivity
function $f\colon\real\to\real$ (with the
$\CIVP$ property or not) which is not almost continuous?}

\pr{pr2}{Does there exist an $f\in\Add\cap\SCIVP\cap\Conn\setminus\AC$?}

\section{An additive almost continuous $\SCIVP$ function
$f\colon\real\to\real$ which is not extendable}

The difficult aspect of constructing a function as in the title will be in
making sure that it will not be extendable.  Since such a function must have
a dense graph (as every discontinuous additive function does) we may
restrict our attention to such functions.  For these we have the following
nice generalization of the $\SCIVP$ property.

\thm{thGenSCIVP}{
If $f\colon\real\to\real$ is an extendable function with a dense graph then
for every $a,b\in\real$, $a<b$, and for each Cantor set $K$ between $f(a)$
and $f(b)$ there is a Cantor set $C$ between $a$ and $b$ such that
$f[C]\subset K$ and the restriction $f\restriction C$ is continuous strictly
increasing.  }

\Proof The basic idea of the proof of this theorem is the same as in the
proof from \cite{RGR} that every extendable function is $\SCIVP$. However,
our schema of the proof will be more similar to the one used to show that
every normal topological space is completely regular.

Let $a$, $b$, and $K$ be as in the theorem and let $\{q_n\colon n<\omega\}$
be an enumeration of some countable subset of $K$ such that the linear
ordering $(\{q_n:n<\omega\},\leq)$ is dense and $q_0=\min K$, $q_1=\max K$.
Since the graph of $f$ is dense (and $f$ is Darboux) we can find $a<b_0<b_1
<b$ with $f(b_0)=q_0$ and $f(b_1)=q_1$.  Let $F\colon\real\times[0,1]\to
\real$ be a connectivity function extending $f$ in a sense that $F(x,0)=
f(x)$ for every $x\in\real$.  By \cite{GR} (see also \cite{64:NatkWoj}) we
can choose $F$ to be continuous outside the line $L_0=\real\times\{0\}$.  We
can also assume that $F(b_0,y)=F(b_0,0)=q_0$ and $F(b_1,y)=F(b_1,0)=q_1$ for
every $y\in[0,1]$. (Indeed, let $H$ be a closed subset of $\real\times[0,1]$
from which we remove two $V$-shape regions with vertices at $\la b_0,0\ra$
and $\la b_1,0\ra$. Extend $F\restriction H$ to $\{0,1\}\times[0,1]$ as
above. Then, by Tietze extension theorem, we can extend such a function to
the reminder of $V$-shape regions continuously. Such modified $F$ will still
be connectivity.)

We will construct a sequence $\la B_n\colon n<\omega\ra$ of compact
connected subsets of $\real\times[0,1]$ such that the following conditions
are satisfied for every $m,n<\omega$,
where $L_\ell=\real\times\{\ell\}$ for
$\ell=0,1$.
\begin{description}
\item{(i)}  $B_0=\{b_0\}\times[0,1]$ and $B_1=\{b_1\}\times[0,1]$.
\item{(ii)} $B_n\cap L_0\neq\emptyset$ and $B_n\cap L_1\neq\emptyset$.
\item{(iii)} If $q_m<q_n$ then, for $\ell=0,1$, we have
\[\max(\{x\in \real:\la x,\ell\ra\in B_m\})<\min(\{x\in\real:\la x,\ell\ra
\in B_n\}).\]
\item{(iv)} $F[B_n]=\{q_n\}$.
\end{description}

Clearly $B_0$ and $B_1$ satisfy (ii)--(iv). So, assume that for some
$n<\omega$, $n>1$, the sets $B_0,\ldots,B_{n-1}$ are already constructed.
To find $B_n$ choose $i,j<n$ such that $(q_i,q_j)$ is the smallest interval
containing $q_n$ with the endpoints from $\{q_0,\ldots,q_{n-1}\}$. Let
\[
b_i=\max(\{x\in \real:\la x,1\ra\in B_i\}),\qquad
  b_j=\min(\{x\in\real:\la x,1\ra\in B_j\}).
\]
(So $b_i<b_j$.) Let $A^*=\cl\big(F^{-1}(q_n)\setminus L_0\big)$. Note that
the set $F^{-1}(q_n)\setminus L_0$ is closed in $\real\times (0,1]$ and thus
$A^*\setminus F^{-1}(q_n)\subseteq L_0$. Now one easily shows that the sets
$B_i$, $B_j$ are contained in different components of the open set
$(\real\times [0,1])\setminus A^*$, so $A^*$ separates $B_i,B_j$. Applying
\cite[thm~4.12, p.\ 51]{Wil} (Property I) we may conclude that there is a
connected component $B^*$ of $A^*$ which separates points
$\la b_i,1\ra$ and $\la b_j,1\ra$,
and thus separates $B_i$ and $B_j$. Note that
$B^*\cap L_0\neq\emptyset\neq B^*\cap L_1$.
Take an $x\in (b_i,b_j)$ such that $\la x,1\ra\in B^*$ and let
$B$ be the connected component of the set $B^*\setminus L_0$ to which $\la
x,1\ra$ belongs. Put $B_n=\cl(B^*)$. We claim that the compact connected set
$B_n$ satisfies our demands. To check clause (iv) note that, by the
definition of the set $A^*$, $B_n\setminus L_0\subseteq F^{-1}(q_n)$. Now
suppose that $y\in L_0\cap B_n$. Assume that $\varepsilon=|F(y)-q_n|>0$.
Since every connectivity function on $\real^2$ is peripherally continuous
(see e.g.\ \cite{GN}), there exists an open neighborhood $W$ of the point $y$
with the diameter $<\frac{1}{2}$ and such that $|F(z)-F(y)|<\varepsilon)$
for all $z\in\bd(W)$.
But $B^*$ is connected, intersects $W\setminus L_0$ and has the
diameter $\geq 1$ ($\cl(B^*)$ intersects $L_0$ and $L_1$), so there exists
a $z\in\bd(W)\cap B^*$, a contradiction. Finally, it should be clear that
$B_n\cap L_1\neq\emptyset$ and $B_n\cap L_0\neq\emptyset$ (e.g.\ use
Proposition~\ref{prop:BBT}), and
\[\begin{array}{l}
\max(\{x\in \real:\la x,\ell\ra\in B_i\})<\min(\{x\in\real:\la x,\ell\ra
\in B_n\}),\\
\max(\{x\in \real:\la x,\ell\ra\in B_n\})<\min(\{x\in\real:\la x,\ell\ra
\in B_j\})
  \end{array}\]
The construction is completed.

Let $B=\bigcup_{n<\omega}B_n$ and notice that
\begin{equation}\label{eqAA}
\mbox{$F\restriction \cl(B)$ is continuous.}
\end{equation}
(Compare \cite[thm~2]{RGR}.) Indeed, by way of contradiction assume that
for some $x\in \cl(B)$ there is a sequence $\la x_i\in B\colon i<\omega\ra$
such that $\lim_{i\to\infty} F(x_i)=L\neq F(x)$. Let $\varepsilon\in(0,
|L-F(x)|/2)$ and $\delta\in(0,1)$ be such that if $|x-x_i|<\delta$ then
$|F(x)-F(x_i)|> \varepsilon$. Using peripheral continuity of the function
$F$ (see e.g. \cite{GN}) we find an open neighborhood $W$ of $x$ with the
diameter $<\delta$ and such that $|f(x)-f(y)|<\varepsilon$ for every
$y\in\bd(W)$. Take $i,n<\omega$ such that $x_i\in W\cap B_n$. Note that
$B_n$ is connected and has the diameter $\geq 1$, so there exists $y\in
\bd(W)\cap B_n$. But then,
\[
\varepsilon<|F(x)-F(x_i)|=|F(x)-y_n|=|F(x)-F(y)|<\varepsilon,
\]
a contradiction.

Consider $L_0$ as ordered in natural order and for $n<\omega$ define
$x_n=\min(B_n\cap L_0)$.  Notice that, by (i) and (iii), $x_n<x_m$ if and
only if $f(x_n)=q_n<q_m=f(x_m)$.  Since $\{q_n\colon n<\omega\}$ (with the
natural order) is a dense linear order, so is $S=\{x_n\colon n<\omega\}$. In
particular, $\cl(S)$ contains a perfect set $C_0=C\times\{0\}$. But $F$ is
continuous on $\cl(B)\supset\cl(S)$ and is strictly increasing on
$S$. Consequently we may choose a perfect set $C^*\subset\cl(S)$ such that
between every two points of $C^*$ there is some $x_n$. So,
$f\restriction C^*$ is
strictly increasing, continuous, and $f[C^*]\subseteq
f[\cl(S)]\subset\cl(\{q_n\colon n<\omega\})=K$. \QED

\thm{thAcScivpNotExt}{There exists
an additive almost continuous $\SCIVP$ function $f\colon\real\to\real$
which is not extendable.}

\proof
Let $\la\la I_\xi,y_\xi\ra\colon\xi<\co\ra$ be a list of all
pairs $\la I,y\ra$ such that $I$ is a nonempty open interval and
$y\in\real$.
Choose the enumerations $\la C_\xi\colon\xi<\co\ra$ of
all perfect subsets of $\real$ and
$\la B_\xi\colon \xi<\co\ra$ of all closed
subsets of $\real^2$ whose projections have non-empty interior.

For our construction we will also use a Hamel basis $H$ which can be
partitioned
onto the sets $\{P_\alpha\colon \alpha\leq\continuum\}$
such that
\begin{itemize}
\item all sets in $\T=\{P_\alpha\colon \alpha<\continuum\}$ are perfect, and

\item
every nonempty open interval contains continuum many $T\in\T$.
\end{itemize}
The existence of such a basis follows easily from the existence
of a linearly independent perfect set
\cite[thm.~2, Ch.~XI sec.~7]{Kucz}
and has been described in details in \cite{SomeDarbF}.

By induction choose a sequence
$\la \la D_\xi,T_\xi\ra\in[H]^{<\omega}\times\T
\colon\xi<\co\ra$
such that the sets
$\{D_\xi\colon \xi<\co\}$
and $\{T_\xi\colon \xi<\co\}$ are pairwise disjoint and that
for every $\xi<\co$
\begin{description}
\item{(i)} $T_\xi\subset I_\xi$,

\item{(ii)} there exists an $a_\xi\in D_\xi\cap\proj(B_\xi)$,

\item{(iii)} there are $z\in \real$,
$0<n<\omega$, non-zero rationals $q_0,\ldots,q_{n-1}$,
and $\{b_0,\ldots,b_{n-1},c_0,\ldots,c_{n-1}\}
\in\left[(D_\xi\setminus\{a_\xi\})\cup\bigcup_{\eta\leq\xi}T_\eta\right]^{2n}$
with the property
that $b_\xi=z+\sum_{j<n}q_j b_j$ and $c_\xi=z+\sum_{j<n}q_j c_j$
belong to $C_\xi$
and that $b_j\in T_\eta$ if and only if $c_j\in T_\eta$ for every
$j<n$ and $\eta\leq\xi$,

\item{(iv)} if $y_\xi\in H$ then
$y_\xi\in\bigcup_{\eta\leq\xi}(D_\eta\cup T_\eta)$.
\end{description}

To make an inductive step assume that
for some $\xi<\co$ the sequence $\la \la D_\eta,T_\eta\ra\colon\eta<\xi\ra$
has been already constructed and let
$M_\xi=\bigcup_{\eta<\xi}(D_\eta\cup T_\eta)$.
It is easy to find $T_\xi\in\T$ with $T_\xi\subset I_\xi\setminus M_\xi$
and an $a_\xi\in\proj(B_\xi)\setminus(T_\xi\cup M_\xi)$.
Next put $\kappa=\left|\bigcup_{\eta<\xi}D_\eta\right|+\omega<\co$ and
for $x\in C_\xi$
let $x=\sum_{i<m_x} q^x_i h^x_i$ be a unique representation
of $x$ in base $H$ (i.e., $q^x_i$'s are non-zero rationals and
$h^x_i$'s are different elements of $H$).
By a combination of the pigeon-hall principle and $\Delta$-system lemma
(see e.g. \cite[thm~1.6, p.~49]{KunenBook})
we can find $m<\omega$,
$\Delta\subset H$, and an $E\subset C_\xi$
of cardinality $\kappa^+$
such that for every different $x,y\in E$ we have:
$m_x=m$, $\Delta=\{h^x_i\colon i<m\}\cap\{h^y_i\colon i<m\}$,
and $q^x_i=q^y_i$ for every $i<m$.
Let $n=m-|\Delta|$.
Refining $E$ and reenumetaring the sets $\{h^x_i\colon i<m\}$, if necessary,
we can also assume that
$h^x_j=h^y_j$ and $\Delta=\{h^x_i\colon n\leq i<m\}$
for all $x,y\in E$ and $n\leq j<m$.
Moreover, since
$|(\xi+2)^n|\leq\kappa<|E|$
we can additionally assume that for every
$i<n$ and $\eta\leq\xi$ we have
$h_i^x\in T_\eta$ if and only if $h_i^y\in T_\eta$.
Finally, by the definition of $\kappa$,
we can also require that
$\{h^x_i\colon i<n\}\cap\left(\{a_\xi\}\cup\bigcup_{\eta<\xi}D_\eta\right)
=\emptyset$ for all $x\in E$.
Fix different $x,y\in E$
and notice that $z=\sum_{n\leq i<m} q^x_i h^x_i$,
$b_i=h^x_i$, $c_i=h^y_i$, and $q_i=q^x_i=q^y_i$ for $i<n$ satisfy (iii).
Now we can
define $D_\xi$ as
$(\{a_\xi\}\cup\{b_0,\ldots,b_{n-1},c_0,\ldots,c_{n-1}\})
\setminus\bigcup_{\eta\leq\xi}T_\eta$
adding to it $y_\xi$, if necessary, to satisfy (iv).
This finishes the inductive construction.

Notice that by (iv) we have
$H=\bigcup_{\xi<\co}(D_\xi\cup T_\xi)$.
We define $f$ on $H$ in such a way
that for each $\xi<\co$ we have:
$\la a_\xi,f(a_\xi)\ra\in B_\xi$, $f\restriction T_\xi\equiv y_\xi$,
and $f(b_i)=f(c_i)$ for every $i<n$, where
$b_i$ and $c_i$ are the points from (iii).
We claim that the unique additive extension of such defined $f\restriction H$
has the desired properties.

Clearly $f$ is additive and almost continuous, since
$f$ intersects every set $B_\xi$.
It is $\SCIVP$ since for every $a<b$
and perfect $K$ between $f(a)$ and $f(b)$
there is $\xi<\co$ with $I_\xi=(a,b)$ and $y_\xi\in K$.
So, $f\restriction T_\xi$ witness $\SCIVP$.
To see that it is not extendable first note that
$f$ is clearly discontinuous, so it has a dense graph.
Thus, by Theorem~\ref{thGenSCIVP} it is enough to show that
$f\restriction C$ is not strictly increasing for every perfect set $C$.
So, let $C$ be perfect.
We claim that there are different $b,c\in C$ such that $f(b)=f(c)$,
which clearly implies that $f\restriction C$ is not strictly increasing.

Indeed let $\xi<\co$ be such that $C=C_\xi$.
Then points
$b_\xi,c_\xi\in C_\xi$ from (iii) are different and
the additivity of $f$ implies that
\[
f(b_\xi)=f\left(z+\sum_{j<n}q_j b_j\right)
=f\left(z+\sum_{j<n}q_j c_j\right)=f(c_\xi).
\]
This finishes the proof. \qed

\section{Another ZFC example of almost continuous $\SCIVP$
function which is not extendable}

In \cite{Ro} H.~Rosen showed that the continuum hypothesis
implies the existence of $\SCIVP$ almost continuous function
$f\colon\real\to\real$ with
a dense graph such that
$f[M]\neq\real$ for every meager set $M\subset\real$.
He also noticed that such an $f$ is not extendable.%
\footnote{In fact, Rosen's function is from $[0,1]$ to $[0,1]$,
but a minor modification gives one from $\real$ to $\real$.}
In this section we will show that a function with such properties
can be constructed in ZFC. (See Theorem~\ref{strangefun}
and Corollary~\ref{AcScivpNotExt}.)
We also show (see Proposition~\ref{bad}) that there
are serious obstacles to make such a function additive.

\lem{dst}{
Suppose that $F\subseteq\baire\times\baire\times\bR$ is a Borel set such that
for some basic open sets $U,V\subseteq\baire$ we have:
\begin{enumerate}
\item[(a)] the set
$Z\stackrel{\rm def}{=}\{\la x,y\ra \in U\times V\colon (F)_{\la x,y\ra }=
\emptyset\}$ is meager,
\item[(b)] the set $A\stackrel{\rm def}{=}\big\{\la x,y\ra \in U\times V\colon
(F)_{\la x,y\ra }\mbox{ is uncountable\/}\big\}$ is meager,
\item[(c)] for each $z\in\bR$ the section $(F)^z$ is meager.
\end{enumerate}
Then there is a perfect set $P\subseteq U\times V$ such that
\[(\forall \la x,y\ra \in P)((F)_{\la x,y\ra }\neq\emptyset)\]
and for distinct $\la x',y'\ra ,\la x'',y''\ra \in P$ we have $(F)_{\la x',y'\ra
}\cap (F)_{\la x'',y''\ra }=\emptyset$ and $x'\neq x''$.
}

\Proof Without loss of generality we may assume that $U=V=\baire$. (Remember
that basic open subsets of $\baire$ are homeomorphic with $\baire$.) Let
$Z^*\subseteq\baire\times\baire$ be a Borel meager set such that $Z\subseteq
Z^*$ and let $A^*\subseteq\baire\times\baire$ be a Borel meager set such
that $A\subseteq A^*$. Take a countable elementary submodel $N$ of $\la
{\cal H}(\chi), \in,<^*\ra$, where $\chi$ is a sufficiently large regular
cardinal, such that the sets $F,Z^*$,
and $A^*$ are in $N$. (Strictly speaking we
require that the Borel codes of these sets are in $N$.)

For $n<\omega$ a set
$T\subseteq\omega^{\textstyle\leq\!n}\times\omega^{\textstyle\leq\!n}$
is an $n$-tree if
\[
\langle\sigma_0,\sigma_1\rangle\in T\ \&\
\sigma_0'\subseteq\sigma_0\
\&\ \sigma_1'\subseteq\sigma_1\quad\Rightarrow\quad \langle\sigma_0',
\sigma_1'\rangle\in T
\]
and for each $\langle\sigma_0,\sigma_1\rangle\in T$ there is $\langle
\sigma_0^*,\sigma^*_1\rangle\in T$ such that
$\sigma_0\subseteq\sigma_0^* \in \omega^{\textstyle n}$ and $\sigma_1
\subseteq\sigma_1^*\in \omega^{\textstyle n}$. Let ${\mathbb P}$ be
the collection of all $T\subseteq\omega^{\textstyle<\!\omega}\times
\omega^{\textstyle<\!\omega}$ which are $n$-trees for some $n$.
We equip
${\mathbb P}$ with the end-extension order,
that is, if $T_0, T_1$ are
$n_0$-- and $n_1$--trees, respectively, then $T_0$ is stronger than $T_1$,
$T_0\leq T_1$, if and only
if $n_1\leq n_0$ and $T_0\cap (\omega^{\textstyle n_1}\times
\omega^{\textstyle n_1})= T_1$. Note that ${\mathbb P}$ is a countable
atomless partial order (and it belongs to $N$), so it is equivalent to the
Cohen forcing notion. Let $G\subseteq{\mathbb P}$ be a generic filter over
$N$.
(It exists since $N$ is countable; of course it is produced by a Cohen real
over $N$.) It is a routine to check that $\bigcup G\subseteq
\omega^{\textstyle <\!\omega}\times \omega^{\textstyle <\!\omega}$ is a
perfect tree. Let
\[P=\{\langle x,y\rangle\in\baire\times\baire: (\forall n\in\omega)(\exists
T\in G)(\langle x\restriction n, y\restriction n \rangle\in T)\}.\]
One easily shows that $P$ is a perfect subset of $\baire\times\baire$ and
that each $\langle x,y\rangle\in P$ is a Cohen real over $N$ (i.e., this
pair does not belong to any meager subset of $\baire\times\baire$ coded in
the model $N$). But even more, all elements of the perfect set $P$ are {\em
mutually\/} Cohen over $N$: if $\langle x,y\rangle$, $\langle x',y'\rangle$
are distinct elements of $P$ then $\langle x,y\rangle$ is Cohen over
$N[\langle x',y'\rangle]$. (Compare with \cite[Lemma 3.3.2]{BaJu95}.)
For our
purposes it is enough to note that if $A\subseteq (\baire\times\baire)^2$ is
a Borel meager set coded in $N$ and $\langle x,y\rangle$, $\langle
x',y'\rangle$ are distinct elements of $P$ then
$\langle\la x,y\ra,\la x',y'\ra\rangle\notin A$.)

We claim that $P$ is as required. First note that $P\cap Z^*=\emptyset$ (so
$(F)_{\la x,y\ra }\neq \emptyset$ for every $\la x,y\ra \in P$) and that
$P\cap A^*=\emptyset$ (implying that $(F)_{\la x,y\ra }$ is countable for
every $\la x,y\ra \in P$). Now suppose that $\la x',y'\ra ,\la x'',y''\ra
\in P$ are distinct.  So $\la x',y'\ra $ is a Cohen real over $N[\la
x'',y''\ra ]$ and in particular $x'\neq x''$ (as $\{x''\}\times\baire$ is a
Borel meager set coded in $N[\langle x'',y''\rangle]$). We know that
$(F)_{\la x'',y''\ra }$ is a countable set from $N[\la x'',y''\ra ]$, and
hence $\bigcup\{(F)^z\colon z\in (F)_{\la x'',y''\ra }\}$ is a meager Borel
set coded in $N[\la x'',y''\ra ]$.  Thus $\la x',y'\ra $ does not belong to
it.  Consequently $(F)_{\la x',y'\ra }\cap (F)_{\la x'',y''\ra }=\emptyset$
and the proof is finished. \QED

\thm{strangefun}{
There is a function $f\colon \bR\to\bR$ such that
\begin{enumerate}
\item[$(\otimes_1)$] if $F\subseteq\bR^2$ is a Borel set such that the
  projection $\proj[F]$ is not meager then $f\cap F\neq\emptyset$,
\item[$(\otimes_2)$] if $P\subseteq\bR$ is a perfect set
  and $B\subseteq\bR$ is a
  non-meager Borel set then there are a perfect set $Q\subseteq B$ and a real
  $y\in P$ such that $f(x)=y$ for all $x\in Q$,
\item[$(\otimes_3)$] if $M\subseteq\bR$ is meager then
  $f[M]\neq\bR$.
\end{enumerate}
}

\Proof First note that $\bR\setminus\bQ$ is homeomorphic to $\baire\times
\baire$, so it is enough to construct a function $f\colon \baire\times\baire
\to\bR$ such that
\begin{enumerate}
\item[$(\otimes_1^*)$] if $F\subseteq(\baire\times\baire)\times\bR$ is a
  Borel set such that the projection of $F$ onto $\baire\times\baire$ is
  not meager then $f\cap F\neq\emptyset$,
\item[$(\otimes_2^*)$] if $P\subseteq\bR$ is a perfect set
  and $B\subseteq\baire\times\baire$
  is a non-meager Borel set then there are a perfect set
  $Q\subseteq B$ and a real $z\in P$ such that $f\rest Q\equiv z$,
\item[$(\otimes_3^*)$] if $M\subseteq\baire\times\baire$ is meager then
$f[M]\not\supset\bR\setminus\{0\}$.
\end{enumerate}
(If a function
$f\colon\bR\setminus\bQ\longrightarrow\bR$ satisfies the demand
$(\otimes_1^*)$--$(\otimes_3^*)$, then the function
$\bar{f}\colon\bR\longrightarrow \bR$ such that $f\subset\bar{f}$ and
$\bar{f}\rest\bQ\equiv 0$ is as required in the theorem.)

Fix enumerations
\begin{itemize}
\item $\{\la r_\alpha,s_\alpha\ra \colon \alpha<\con\}$ of $\baire\times\baire$,
\item $\{M_\alpha\colon \alpha<\con\}$ of all Borel meager subsets of
$\baire\times
\baire$,
\item $\{\la P_\alpha,B_\alpha\ra \colon \alpha<\con\}$ of pairs $\la
P,B\ra $ such
that
$P\subseteq\bR$ is a perfect set, and $B\subseteq\baire\times\baire$ is a
Borel non-meager set,
\item $\{F_\alpha\colon \alpha<\con\}$ of all Borel sets $F\subseteq\baire\times
\baire\times\bR$ such that the projection of $F$ onto $\baire\times\baire$ is
not meager.
\end{itemize}
By induction on $\alpha<\con$ we will choose perfect sets $Q_\alpha\subseteq
\baire$ and reals $x^0_\alpha,x^1_\alpha,w_\alpha\in\baire$, $y_\alpha,
z_\alpha\in\bR$,
$v_\alpha\in\bR\setminus\{0\}$ such that for $\alpha,\beta<\con$:
\begin{enumerate}
\item[(i)]   $(\{w_\alpha\}\times Q_\alpha)\subseteq B_\alpha$, $z_\alpha\in
P_\alpha$,
\item[(ii)]  $\la x_\alpha^0,x_\alpha^1\ra \notin \{w_\beta\}\times Q_\beta$
and if $\alpha\neq\beta$ then $\la x_\alpha^0,x_\alpha^1\ra \neq \la x_\beta^0,
x_\beta^1\ra $, $w_\alpha\neq w_\beta$, and $v_\alpha\neq v_\beta$,
\item[(iii)] $\la x_{2\alpha}^0,x_{2\alpha}^1,y_{2\alpha}\ra \in F_\alpha$,
\item[(iv)]  $\la r_\alpha,s_\alpha\ra \in\{\la x_\gamma^0,x_\gamma^1\ra \colon
\gamma\leq 2\alpha+1\}\cup\bigcup\limits_{\gamma\leq 2\alpha+1}
\{w_\gamma\}\times
Q_\gamma$,
\item[(v)]   $z_\alpha\neq v_\beta$, and if $\la x_\alpha^0,x_\alpha^1\ra \in
M_\beta$ then $y_\alpha\neq v_\beta$.
\end{enumerate}
Assume that we can carry out the construction so that the demands (i)--(v) are
satisfied. Define a function $f\colon \baire\times\baire\to\bR$ by:
\[(\forall\alpha<\con)(f\rest (\{w_\alpha\}\times Q_\alpha)\equiv z_\alpha\
\&\ f(x_\alpha^0,x_\alpha^1)= y_\alpha).\]
It follows from the clauses (ii) and (iv), that the above condition defines a
function on $\baire\times\baire$. This function has the required properties:
$(\otimes_1^*)$ holds by clause (iii), $(\otimes_2^*)$ follows from clause (i),
and $(\otimes_3^*)$ is a consequence of (v) since $v_\alpha\notin
f[M_\alpha]$.

So let us show how the construction may be carried out. Assume that we have
defined $x_\beta^0, x_\beta^1,w_\beta\in\baire$, $y_\beta,z_\beta\in\bR$,
$v_\beta\in\bR\setminus\{0\}$, and $Q_\beta\subseteq\baire$
for $\beta<\alpha$.
First choose non-zero numbers
$v_\alpha\in
\bR\setminus\bigcup\{\{v_\beta,y_\beta,z_\beta\}\colon\beta<\alpha\}$
and $z_\alpha\in P_\alpha\setminus\{v_\beta\colon \beta\leq\alpha\}$.
The set $B_\alpha$ is not meager so we find
$w_\alpha\in\baire
\setminus\bigcup\{\{x_\beta^0,w_\beta\}\colon \beta<\alpha\}$
such
that the section
$(B_\alpha)_{w_\alpha}$ is not meager. Pick a perfect set $Q_\alpha\subseteq
(B_\alpha)_{w_\alpha}$. Next we consider two separate cases to choose
$x_\alpha^0$, $x_\alpha^1$, and $y_\alpha$.

\medskip

\noindent{\sc Case 1:} $\alpha$ is odd, say $\alpha=2\alpha_0+1$.

\medskip

Let $\la x_\alpha^0,x_\alpha^1\ra \in\baire\times\baire\setminus
\left(\{\la x_\beta^0,
x_\beta^1\ra \colon \beta<\alpha\}\cup\bigcup\limits_{\beta\leq\alpha}
\{w_\beta\}
\times Q_\beta\right)$
be such that
\[\la r_{\alpha_0},s_{\alpha_0}\ra \in\{\la x_\beta^0,x_\beta^1\ra \colon
\beta\leq\alpha\}\cup
\bigcup\limits_{\beta\leq\alpha} \{w_\beta\}\times Q_\beta,\]
and let $y_\alpha\in\bR\setminus\{v_\beta,z_\beta\colon \beta\leq\alpha\}$.

\medskip

\noindent{\sc Case 2:} $\alpha$ is even, say $\alpha=2\alpha_0$.

\medskip

Look at the set $F_{\alpha_0}$. If there is $y\in\bR$ such that the section
$\big(F_{\alpha_0}\big)^y$ is not meager then take such an $y$ as $y_\alpha$.
Pick
\[
x_\alpha^0\in\baire\setminus
\bigcup\{\{w_\beta,x_\beta^0\}\colon \beta<\alpha\}\setminus\{
w_\alpha\}
\]
such that $\big(\big(F_{\alpha_0}\big)^{y_\alpha}\big)_{x_\alpha^0}$ is not
meager and
\begin{quotation}
\noindent if $v_\beta=y_\alpha$, $\beta\leq\alpha$, then
$(M_\beta)_{x_\alpha^0}$ is meager.
\end{quotation}
(Note that there is at most one $\beta$ as above.) Next choose $x_\alpha^1\in
\baire$ such that $\la x_\alpha^0,x_\alpha^1,y_\alpha\ra \in F_{\alpha_0}$
and
$\la x_\alpha^0,x_\alpha^1\ra \notin M_\beta$ provided $v_\beta=y_\alpha$,
$\beta\leq\alpha$.

So suppose now that for each $y\in\bR$ the section $\big(F_{\alpha_0}\big)^y$
is meager. Let
\[
A\stackrel{\rm def}{=}\{\la x_0,x_1\ra \in\baire\times\baire\colon
(F_{\alpha_0})_{\la  x_0,x_1\ra }\mbox{ is uncountable\/}\}.
\]
It is an analytic set,
so it has the Baire property. If $A$ is not meager
then we may choose
$x_\alpha^0\in \baire\setminus\bigcup\{\{w_\beta,x_\beta^0\}\colon \beta<
\alpha\}\setminus\{w_\alpha\}$
and $x_\alpha^1\in\baire$ and $y_\alpha\in\bR
\setminus\{v_\beta\colon \beta\leq\alpha\}$ such that $\la
x_\alpha^0,x_\alpha^1,
y_\alpha\ra \in F_{\alpha_0}$.

So assume that the set $A$ is meager. Take basic open sets $U,V\subseteq
\baire$ such that $\{\la x_0,x_1\ra \in U\times V\colon  (F_{\alpha_0})_{\la
x_0,x_1\ra }=
\emptyset\}$ is meager. Note that the sets $U,V$ and $F_{\alpha_0}$ satisfy
the assumptions of
Lemma~\ref{dst}. So we get a perfect set $P\subseteq U\times V$
such that $(F_{\alpha_0})_{\la x_0,x_1\ra }\neq\emptyset$
for every $\la x_0,x_1\ra \in P$ and that for distinct
$\la x_0',x_1'\ra ,\la x''_0,x''_1\ra\in P$:
\[
(F_{\alpha_0})_{\la x_0',x_1'\ra }\cap (F_{\alpha_0})_{\la
x''_0,x''_1\ra }=\emptyset\quad\mbox{ and }\quad x'_0\neq x''_0.
\]
Now we may easily find $\la x_\alpha^0,x_\alpha^1\ra \in P$ and $y_\alpha\in\bR
\setminus\{v_\beta\colon \beta\leq\alpha\}$ such that
\[x_\alpha^0\notin\{w_\beta,x_\beta^0\colon \beta<\alpha\}\cup\{w_\alpha\}\quad
\mbox{ and }\quad \la x_\alpha^0,x_\alpha^1,y_\alpha\ra \in F_{\alpha_0}.\]
This finishes the inductive step of the construction. Checking that the
demands (i)--(v) are satisfied is straightforward in all cases. (Note that
it follows from $(\otimes^*_1) + (\otimes^*_3)$ that for each meager set
$M\subset\baire\times\baire$, the set $\bR\setminus f[M]$ is
uncountable. One may easily guarantee that these sets are of size $\con$,
but there is no need for this.)

Thus the proof of the theorem is complete. \QED

\cor{AcScivpNotExt}{
There exists an almost continuous function $f\colon \bR \to\bR$ which
has the strong Cantor intermediate value property but is not an extendability
function.
}

\Proof Let $f\colon \bR\to\bR$ be the function constructed in
Theorem~\ref{strangefun}.
The property $(\otimes_1)$ implies that the function $f$ is
almost continuous and $\rng(f)$ is dense in $\bR^2$, and the property
$(\otimes_2)$  guarantees that $f\in\SCIVP$. To show that $f$ is not an
extendability function we use the third property listed in
Theorem~\ref{strangefun}.
So by way of contradiction assume
 that $f\in\ext$. Then, by Rosen \cite{Ro94}, there is a meager set
$M\subseteq\bR$ such that
\begin{enumerate}
\item[$(\oplus)$] if $g\colon \bR\to\bR$ and $g\rest M=f\rest M$ then $g$
is an extendability function.
\end{enumerate}
We may additionally demand that $\cl(f[M])=\bR$. (Just increase $M$ if
necessary.)
Pick any $r^*\in f[M]$ and define a function $g\colon\bR\to\bR$ by:
\[
g(x)=\left\{\begin{array}{ll}
f(x)&\mbox{ if }x\in M,\\
r^* &\mbox{ otherwise.}
              \end{array}\right.
\]
By $(\oplus)$, $g$ is an extendability function (and thus Darboux) and by its
definition $\rng(g)=f[M]$ is a dense subset of $\bR$ (so it has to be $\bR$).
But $f[M]\neq\bR$ (remember $(\otimes_3)$ of
Theorem~\ref{strangefun}), a
contradiction.  \QED
\medskip

One would hope for getting an additive function as in
Theorem~\ref{strangefun}. Unfortunately this approach cannot work.

\begin{proposition}
\label{bad}
Suppose that $f\colon \bR\to\bR$ is an additive function such that
\begin{enumerate}
\item for some perfect set $P\subseteq\bR$, the restriction $f\rest P$ is
continuous,
\item for each nowhere dense set $M\subseteq\bR$, the image $f[M]$ is not
$\bR$.
\end{enumerate}
Then there is a closed set $F\subseteq\bR^2$ such that $\proj[F]=\bR$ and
$f\cap F=\emptyset$.
\end{proposition}

\Proof Let $P\subseteq\bR$ be a compact perfect set such that $f\rest P$ is
continuous. By Erd\H os, Kunen, and Mauldin \cite{EKM}, we find a compact
perfect set $Q$ of Lebesgue measure 0 (and so nowhere dense) such that
$P+Q$ contains the interval $[0,1]$. By the second assumption, we may pick a
real $r\in \bR\setminus f[Q]$. Let $F$ be the subset of the plane $\bR^2$
described by:
\[\begin{array}{l}
\la x,y\ra\in F\qquad\mbox{if and only if}\\
(\exists w\in P)(\exists z\in Q)(\exists m\in\bZ)(x=w+z+m\ \&\ y=f(w)+f(m)+r).
  \end{array}\]
Since $P,Q$ are compact and $f\rest P$ is continuous, the set $F$ is closed.
By the choice of the perfect $Q$ we know that $\proj[F]=\bR$. Finally, suppose
that $\la x,y\ra\in F\cap f$. Take $w\in P$, $z\in Q$ and $m\in\bZ$ witnessing
$\la x,y\ra\in F$. Then
\[f(w)+f(m)+r=y=f(x)=f(w+z+m)=f(w)+f(z)+f(m),\]
and hence $f(z)=r$, a contradiction with the choice of $r$. \QED


\begin{thebibliography}{99}
\bibitem{BCN} M.~Balcerzak, K.~Ciesielski, T.~Natkaniec,
\newblock{\em Sierpi\'nski--Zygmund functions that are Darboux, almost
continuous, or have a perfect road},
\newblock Arch. Math. Logic {\bf 31}(1) (1998), 29--35.
(Preprint$^\star$ available.%
\footnote{Preprints marked by $^\star$ are available in an electronic
form. They can be accessed from Set Theoretic Analysis web page:\\
http://www.math.wvu.edu/homepages/kcies/STA/STA.html})

\bibitem{Ban} D.~Banaszewski,
\newblock {\em On some subclasses of additive functions},
\newblock PhD Thesis, {\L}\'od\'z University 1997 (in Polish).

\bibitem{KBanNat} K.~Banaszewski, T.~Natkaniec, {\it Sierpi\'nski--Zygmund
functions that have the Cantor intermediate value property},
Real Anal. Exchange, to appear.
(Preprint$^\star$ available.)

\bibitem{BaJu95} T.~Bartoszy\'nski, H.~Judah,
\newblock {\em {Set Theory: On the Structure of the Real Line}}.
\newblock A K Peters, Wellesley, Massachusetts 1995.

\bibitem{ChVo} C.~O.~Christenson, W.~L.~Voxman,
\newblock {\em Aspects of topology}, volume 39 of Monographs and Textbooks in
Pure and Applied Mathematics,
\newblock Marcel Dekker Inc., New York -- Basel 1977.

\bibitem{Ci98} K.~Ciesielski,
\newblock {\em Set Theoretic Real Analysis},
\newblock {J. Appl. Anal.} {\bf 3}(2) (1997), 143--190.
(Preprint$^\star$ available.)

\bibitem{CiBook} K.~Ciesielski,
\newblock {\em Set Theory for the Working Mathematician},
\newblock London Math. Soc. Student Texts {\bf 39}, Cambridge Univ. Press
1997.

\bibitem{SomeDarbF}
{\em Some additive Darboux-like functions},
{J. Appl. Anal.} {\bf 4}(1) (1998), 43--51.
(Preprint$^\star$ available.)


\bibitem{CiJa} K.~Ciesielski, J.~Jastrz{\c{e}}bski,
\newblock {\em Darboux--like functions within the classes of Baire one, Baire
two, and additive functions}, Topology Appl., to appear.
(Preprint$^\star$ available.)

\bibitem{64:NatkWoj} K.~Ciesielski, T.~Natkaniec, J.~ Wojciechowski,
{\em  Extending connectivity functions on $\real^n$}, preprint$^\star$.


\bibitem{EKM} P.~Erd\H os, K.~Kunen, R.~D.~Mauldin,
\newblock {\em Some additive properties of sets of real numbers},
\newblock Fund. Math. {\bf 113} (1981), 187--199.

\bibitem{GN} R.~G.~Gibson, T.~Natkaniec,
\newblock {\em Darboux like functions},
\newblock Real Anal. Exchange {\bf 22}(2) (1996--97),
492--533. (Preprint$^\star$ available.)

\bibitem{GR} R.~G.~Gibson, F.~Roush, {\em A characterization
of extendable connectivity functions}, Real Anal. Exchange
{\bf 13} (1987--88), 214--222.

\bibitem{ZG} Z.~Grande,
\newblock {\em On almost continuous additive functions},
\newblock Math. Slovaca, to appear.

\bibitem{GMN}  Z.~Grande, A.~Maliszewski, T.~Natkaniec,
\newblock {\em Some problems concerning almost continuous functions},
\newblock  Proceedings of the Joint US-Polish Workshop in Real Analysis, Real
Anal. Exchange {\bf 20} (1994--95), 429--432.

\bibitem{Jo} F.~Jordan,
\newblock {\em Cardinal invariants connected with adding real functions}.
\newblock {Real Anal.  Exchange} {\bf 22} (1996--97), 696--713.
(Preprint$^\star$ available.)

\bibitem{Kucz} M. Kuczma,
\newblock {\em An Introduction to the Theory of Functional Equations and
Inequalities},
\newblock Polish Scientific Publishers PWN, Warsaw 1985.

\bibitem{KunenBook}  K.~Kunen, {\it Set Theory}, North Holland, Amsterdam 1980.

\bibitem{Kur} K. Kuratowski,
\newblock {\em Topology},
\newblock Vol. II, Acad. Press, New York, N.Y., 1968.

\bibitem{Nadler} S.~B.~Nadler, Jr.,
\newblock {\em Continuum Theory},
\newblock Marcel Dekker, Inc., New York, N.Y., 1992.

\bibitem{N1} T.~Natkaniec,
\newblock {\em Almost continuity},
\newblock Real Anal. Exchange {\bf 17} (1991--92), 462--520.

\bibitem{Ro65} J.~H.~Roberts,
\newblock {\em Zero--dimensional sets blocking connectivity functions}.
\newblock  Fund. Math. {\bf 57} (1965), 173--179.

\bibitem{Ro94} H.~Rosen,
\newblock {\em Limits and sums of extendable connectivity functions},
\newblock {Real Anal. Exchange} {\bf 20} (1994--95), 183--191.

\bibitem{Ro} H.~Rosen,
{\em An almost continuous nonextendable function},
Real Anal. Exchange {\bf 23}(2) (1997--98), 567--570.

\bibitem{RGR}
H.~Rosen, R.~G.~Gibson, F.~Roush, {\em Extendable functions
and almost continuous functions with a perfect road}, Real Anal.
Exchange {\bf 17} (1991--92), 248--257.

\bibitem{Stal} J.~Stallings,
\newblock {\em Fixed point theorems for connectivity maps}.
\newblock Fund. Math. {\bf 47} (1959), 249--263.

\bibitem{Wil} R.~L.~Wilder, {\em Topology of Manifolds},
AMS Colloquium Publ. {\bf 32} (1949).

\end{thebibliography}
\end{document}